\documentclass[10pt, conference, final]{ieeeconf}
\IEEEoverridecommandlockouts

\usepackage{lmodern}
\usepackage[T1]{fontenc}
\usepackage[utf8]{inputenc}							
\usepackage{amsmath}									
\usepackage{amssymb}
\usepackage[usenames,dvipsnames]{color}
\usepackage{multicol}
\usepackage[lined,boxed]{algorithm2e}
\usepackage[hidelinks]{hyperref}
\usepackage{graphicx}	
\let\labelindent\relax
\usepackage{enumitem}
\usepackage{caption}
\usepackage{subcaption}

\usepackage{epstopdf} 
	\DeclareGraphicsExtensions{.pdf,.png,.jpg,.eps}

\usepackage{cite}
\usepackage{xcolor}

\newcommand{\mset}{\mathbf{\Sigma}}
\newcommand{\graph}{\Theta}
\newcommand{\reels}{\mathbb{R}}

\newcommand{\jsr}{\hat \rho}

\newtheorem{theorem}{ Theorem}[section]
\newtheorem{conjecture}{ Conjecture}
\newtheorem{remark}{Remark}
\newtheorem{assumption}{Assumption}
\newtheorem{definition}{Definition}

\newtheorem{proposition}[theorem]{Proposition}

\newcommand{\s}[1]{\underset{#1}{\sup}}
\newcommand{\im}{\text{Im}}
\renewcommand{\ker}{\text{Ker}}

\newcommand{\Ell}{\mathcal{L}}

\usepackage{comment}

\newif\ifshortVersion
\shortVersionfalse 
\ifshortVersion
	\excludecomment{shortVersion}
\else
	\newenvironment{shortVersion}{\sffamily\color{blue}}{}
\fi

\date{}
\title{Extremal storage functions and minimal realizations of discrete-time linear switching systems}
 
\author{Matthew Philippe,  Ray Essick,  Geir Dullerud and Rapha\"{e}l M. Jungers%
\thanks{M. Philippe is a FRIA (F.R.S. - FNRS) Fellow; R. Jungers is a FNRS Research Associate. Both are with the ICTEAM institute from the Universit\'{e} catholique de Louvain. Louvain-la-Neuve, 1348, Belgium. \{matthew.philippe, raphael.jungers\}@uclouvain.be}%
\thanks{R. Essick and G. E. Dullerud were partially supported by grants NSA SoS W911NSF-13-0086 and AFOSR MURI FA9550-10-1-0573. Both are with the MecSE and CSL, University of Illinois at Urbana-Champaign. Urbana, IL 61801, USA. \{ressick2, dullerud\}@illinois.edu.}%
}
\begin{document}
\maketitle
\begin{abstract}
We study the $\Ell_p$ induced gain of discrete-time linear switching
systems with graph-constrained switching sequences.
 We first prove that, for stable systems in a minimal realization, for
 every $p \geq 1$, the $\Ell_p$-gain is exactly
 characterized through switching storage functions.
 These functions are shown to be the $p$th power of a norm.
 In order to consider general systems, we provide an algorithm for computing minimal realizations.
  These realizations are  \emph{rectangular systems}, with a state
  dimension that varies according to the mode of the system.\\
   We apply our tools to the study on the of $\Ell_2$-gain. We provide algorithms for its approximation, and provide a converse result for the existence
   of quadratic switching storage functions. We finally illustrate the results with a  physically motivated example.
\end{abstract}
\section{Introduction}
Discrete-time linear switching systems are multi-modal systems that switch between a finite set of modes. They arise in many practical and theoretical applications \cite{LiSISA,JuTJSR,AlRaCMAA,JuDiFSOD,HeMiOAMS,ShWiAPSM}.
They are of the form
 \begin{equation}
\begin{aligned}
x_{t+1} & = A_{\sigma(t)} x_t + B_{\sigma(t)} w_t, \\
 z_{t} & = C_{\sigma(t)} x_t + D_{\sigma(t)} w_t,
\end{aligned}
\label{eq:switching}
\end{equation} 
where $x \in \reels^n$, $w \in \reels^d$ and $z \in \reels^m$ are respectively the state of the system, 
a \emph{disturbance input}, and a \emph{performance output}. The function $\sigma(\cdot) : \mathbb{N} \rightarrow \{1, \ldots, N\}$ is called a \emph{switching sequence}, and at time $t$, $\sigma(t)$ is called the \emph{mode} of the system at time $t$. We consider  switching sequences  following a set of logical rules: they need be generated by walks in a given labeled graph $\graph$ (we will denote this by $\sigma \in \graph$).\\
The analysis and control of switched systems is an active research area, and many questions that are easy to understand and decide for LTI (i.e., Linear Time Invariant) systems are known to be hard for switching systems (see, e.g., \cite{BlTsTBOA, JuTJSR, LiSISA, PhMiDTBA}). Nevertheless, several analysis and control techniques have been devised (see e.g. \cite{LeDuUSOD, LeDuODAF, KuChSSSF, EsLeCOLS}), often providing conservative yet tractable methods. The particular question of the \emph{stability} of a switched system attracts a lot of attention; in the past few years, approaches using multiple Lyapunov functions have been devised \cite{DaBePDLF, BlFeSAOD,LeDuUSOD}, with recent works \cite{PhEsSODT, EsPhTMAS, AhJuJSRA} analyzing how conservative these methods are.
\\
In this paper, we are interested in the analysis and computation of the \emph{$\Ell_p$ induced} gain of switching systems:
\begin{equation}
\gamma_p = \left ( \s{ \mathbf{w} \in \Ell_p, \, \sigma \in \graph, \, x_0 = 0.} \frac{\sum_{t = 0}^{\infty} \|z_t\|^p_p}{ \sum_{t = 0}
^{\infty} \|w_t\|^p_p } \right )^{1/p},
\label{eq:gain}
\end{equation}
where $\mathbf{w} \in \Ell_p$ means that the disturbance signal satisfies $\sum_{t = 0}^{\infty} \|w_t\|_p^p < + \infty$, $\|w\|_p$ being the $p$-norm of $w$.\\ This quantity is a  measure of a system's ability to dissipate the energy of a disturbance.
As noted in \cite{NaVoAUFF}, the $\Ell_2$-gain is one of the most studied performance measures. Some approaches link the gain of a switching system to the individual performances of its LTI  modes (see e.g. \cite{CoBoRMSG, HeLIGO}); in \cite{EsLeCOLS, LeDuODAF}, a  hierarchy of LMIs is presented to decide whether $\gamma_2 \leq 1$. In \cite{PuZhASOT}, a generalized version of the gain is studied through generating functions. For continuous time systems, \cite{HiHeLIGA} gives a tractable approach based on the computation of a storage function. For the  general \emph{$\Ell_p$-gain}, we refer to \cite{NaVoCAOO} that relies on an operator-theoretic approach.\\
Our approach is inspired from the works \cite{PhEsSODT, EsPhTMAS, AhJuJSRA} on \emph{stability analysis}. These works present a framework for the stability analysis of switching systems (with constrained switching sequences) using multiple Lyapunov functions with pieces that are \emph{norms}. In Section \ref{sec:StorageTheory}, we provide a characterization of $\Ell_p$-gain using \emph{multiple storage functions}, where pieces are now the $p$th power of norms. We rely on two assumptions, namely that
the switching system be \emph{internally stable} and in a \emph{minimal realization}. In order to assert the generality of the results, both assumptions are discussed  in Subsections \ref{subsec:Stability} and Subsection \ref{subsec:Minimize} respectively. In particular, we give an algorithm for computing minimal realizations: these are  \emph{rectangular systems} for which the dimension of the state space varies with the modes of the system. In Section \ref{sec:L2}, we narrow our focus to the $\Ell_2$-gain, providing two approaches for estimating the gain and constructing storage functions. The first uses dynamic programming and provides asymptotically tight lower bound on the gain, while the second, based on the work of \cite{EsLeCOLS, LeDuODAF}, provides  asymptotically tight upper bounds.
Then, in Subsection \ref{Subsection:Converse} we present a converse result for the existence of \emph{quadratic storage functions}, which can be  obtained using convex optimization.
 Section \ref{sec:Example}  illustrates our results with a simple practical example.
\section*{Preliminaries}
\begin{figure}[!ht]
\centering
\begin{subfigure}[b]{0.2\textwidth}
\centering
\includegraphics[scale = 0.4]{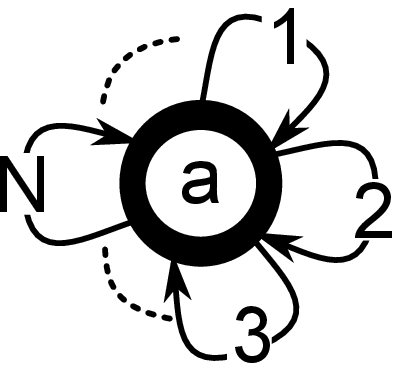}
\caption{}
\label{fig:intro_exb}
\end{subfigure}
~
\begin{subfigure}[b]{0.2\textwidth}
\centering
\includegraphics[scale = 0.4]{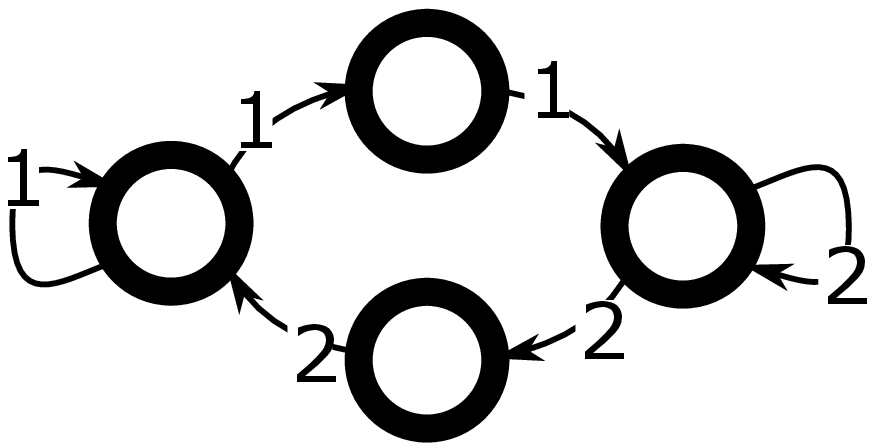}
\caption{}
\label{fig:intro_exa}
\end{subfigure}
\caption{ The automaton of Fig. \ref{fig:intro_exb} generates any switching sequences on $N$ modes. Fig. \ref{fig:intro_exa} generates the sequence $\{1,1,2,2,1,1, \ldots, \}$, but cannot generate $\{1,2,1,2, \ldots \}$.  }
\label{fig:intro_ex}
\end{figure}

Given a system of the form (\ref{eq:switching}) on $N$ modes, we denote the set of parameters of the system by
$$ \mset= \{ (A_\sigma, B_\sigma, C_\sigma, D_\sigma)_{\sigma = 1, \ldots, N}\}. $$  
We define constraints on the switching sequences through a strongly connected, labeled and directed graph $\graph(V,E)$ (see e.g.  
\cite{PhEsSODT, KuChSSSF, EsLeCOLS}).  The edges of $\graph$ are of the form $(u,v,\sigma) \in E$, where $u,v \in V$ are the source and destination nodes and $\sigma \in \{1, \ldots, N\}$ is a label that corresponds to a mode of the switching system. A \emph{constrained switching system} is then defined by a pair $(\graph, \mset)$, and admissible switching sequences correspond to paths in $\graph$ (see Fig. \ref{fig:intro_ex} for examples).
More precisely, let $\pi \in \graph$ denote a finite-length path $\pi$ in the graph. The associated switching sequence is written $t \mapsto \sigma_{\pi}(t)$, where, for a given $t \geq 0$, $\sigma_{\pi}(t)$ is the \emph{label} on the $t+1$th edge of $\pi$.  The length of $\pi$, i.e., the number of \emph{edges} it contains, is written $|\pi|$. We let $\pi(i:j)$ with $1 \leq i \leq j \leq |\pi|$ denote the path formed from the $i$th edge of $\pi$ to its $j$th (both included).  
We let $t \mapsto v_{\pi}(t) \in V$ be the sequence of \emph{nodes} encountered by $\pi$, where $v_{\pi}(0)$ is the source node, and if $|\pi| < + \infty$, $v_{\pi}(|\pi|)$ is the destination node. The functions $v_\pi(t)$ and $\sigma_{\pi}(t)$ need not be defined for $t \geq |\pi|$. In order to compactly describe the input-output maps of the system, we write
\begin{equation}
A_\pi := A_{\sigma_\pi(|\pi|-1)} \cdots A_{\sigma_\pi(0)},
\label{eq:Amap}
\end{equation}
\begin{equation}
B_\pi := \begin{pmatrix} A_{\sigma_\pi(|\pi|-1)}\cdots A_{\sigma_\pi(1)}B_{\sigma_\pi(0)},& \ldots,  &   B_{\sigma_\pi(|\pi|-1)}\end{pmatrix},
\label{eq:Bmap}
\end{equation}
\begin{equation}
C_\pi := \begin{pmatrix}
    C_{\sigma_\pi(0)},\\
    C_{\sigma_\pi(1)}A_{\sigma_\pi(0)},\\
    \vdots \\
    C_{\sigma_\pi(|\pi|-1)}A_{\sigma_\pi(|\pi|-2)} \cdots A_{\sigma_\pi(0)} 
\end{pmatrix},
\label{eq:Cmap}
\end{equation}
\begin{multline}
D_\pi := \\
\begin{pmatrix}
   D_{\sigma_\pi(0)}                   & \hspace{-12pt} 0          & \hspace{-9pt} & \hspace{-6pt} 0 \\
   C_{\sigma_\pi(1)}B_{\sigma_\pi(0)} & \hspace{-12pt} D_{\sigma(1)} & \hspace{-9pt} & \hspace{-6pt} 0 \\
   \vdots                                & \hspace{-12pt} \vdots        & \hspace{-9pt} \ddots & \hspace{-6pt} 0 \\
   C_{\sigma_\pi(|\pi| - 1 )}A_{\sigma_\pi(|\pi| - 2 )} \cdots A_{\sigma_\pi(1)}B_{\sigma_\pi(0)} & \hspace{-12pt} \ldots & \hspace{-9pt} \ldots & \hspace{-6pt} D_{\sigma_\pi(|\pi| - 1)}
\end{pmatrix}\hspace{-1pt}\\.
\label{eq:Dmap}
\end{multline}
If $\pi$ is an \emph{infinite} path, we define $A_\pi$, $B_\pi$, $C_\pi$, and $D_\pi$ similarly, the last three being operators on $\Ell_p$.
When in a matrix, $0$ and $I$ are the null and identity matrix of appropriate dimensions. \\
Given a matrix $M \in \reels^{n \times m}$, $\im(M)$ is the image of $M$ and $\ker(M)$ the kernel of $M$. The orthogonal of a subspace $\mathcal{X}$ of $\reels^n$ is $\mathcal{X}^\perp$. The sum of $K$ subspaces of $ \reels^n$, $(\mathcal{X}_i)_{i \in 1, \ldots, K}$,  is $\sum_{i=1}^K \mathcal{X}_i = \{ \sum_{i=1}^K x_i: \, x_i \in \mathcal{X}_i\}.$ 

\section{Extremal storage functions}
\label{sec:StorageTheory}
Our main results for this section rely on two assumptions. 
They are discussed in Subsections \ref{subsec:Stability} and \ref{subsec:Minimize}. 
\begin{assumption}[Minimality]
The system $(\graph(V,E), \mset)$ is \emph{minimal} in the sense that, defining for all $v \in V$,
\begin{equation}
\mathcal{B}_v = \hspace{-6pt} \sum_{\pi \in \graph, \, v_\pi(|\pi|) = v} \hspace{-6pt} \im(B_\pi), \ \ \mathcal{C}_v = \hspace{-6pt} \bigcap_{\pi \in \graph, \, v_\pi(0) = v} \hspace{-6pt}
\ker(C_\pi),
\label{ass:eq:reach}
\end{equation} where $B_\pi$ and $C_\pi$ are given in (\ref{eq:Bmap}), (\ref{eq:Cmap}), we have 
$\mathcal{B}_v = \reels^{n}, \, \mathcal{C}_v = \{0\} \subset \reels^{n}. $
\label{ass:minimal}
\end{assumption}
 This definition of minimality equivalent to that found in \cite{PeBaRTOD}, Theorem 1.
Minimality requires that any state $x \in \reels^n$ has a reachable and a detectable component. 
\begin{assumption}[Internal stability]
The system $(\graph, \mset )$ is \emph{internally (exponentially) stable}, i.e., 
there exist $K \geq 1$ and $\rho < 1$ such that $ \forall \pi \in \graph, \|A_\pi\| \leq K \rho^{|\pi|}.$
\label{ass:stability}
\end{assumption}
Internal stability guarantees that the $\Ell_p$-gain of the system is bounded (see e.g. \cite{DuLaANAF}). It is difficult to decide whether a given constrained system is internally stable. Towards this end,
generalizing the tools introduced in  \cite{AhJuJSRA}, 
the concept of \emph{multinorm} is introduced in \cite{PhEsSODT}:
\begin{definition}
A multinorm for a system $(\graph(V,E), \mset)$ is a set   $\{(|\cdot|_v)_{v \in V}\}$ of one norm of $\reels^n$ per node in $V$. 
\end{definition} 
This concept allowed to develop tools for the approximation of the \emph{constrained joint spectral radius}, denoted $\jsr$,  of systems $(\graph, \mset)$ (see \cite{DaAGTS,PhEsSODT}).
It can be shown that the internal stability of a system is equivalent to $\hat{\rho} < 1$ (see e.g. \cite{DaAGTS,KoTBWF}), and in \cite{PhEsSODT}-Proposition 2.2, the authors show that
\begin{equation}
\jsr = \inf \rho : \left \{
\begin{aligned}
&\exists \{(|\cdot|_v)_{v \in V}\} : \\
& \forall x \in \reels^n,  \forall (u,v,\sigma) \in E, |A_\sigma x |_v  \leq \rho|x|_u.
\end{aligned} \right .
\label{eq:jsrmulno}
\end{equation}
Thus, a constrained switching system is stable if and only if it has a multiple Lyapunov function, with no more than one piece per node of $\graph$, which can be taken as a norm. Whenever the infimum is actually attained by a multinorm, we say it is \emph{extremal}. The existence of such an extremal multinorm is not guaranteed. 

We now show that multinorms provide a characterization of the $\Ell_p$-gain of constrained switching systems. They  play the role of \emph{multiple storage functions}.
\begin{theorem}
Consider a system $(\graph, \mset)$ satisfying Assumptions \ref{ass:minimal} and  \ref{ass:stability}. Its $\mathcal{L}_p$-gain satisfies
\begin{align}
\gamma_p(\graph, & \mset) = \min \gamma: \label{eq:optigain} \\ 
         &  \exists \{(|\cdot|_v)_{ v \in V}\}\ s.t.\ \forall x \in \reels^n, \forall w, \forall (u,v,\sigma) \in E,  \nonumber \\
         & 
         \left ( 
         \begin{aligned}
          |A_\sigma x + B_\sigma w|_v^p + \|C_\sigma x + D_\sigma w\|_p^p   \\
          \leq |x|_u^p + \gamma^p \|w\|_p^p. 
         \end{aligned} \label{eq:normineqgain}
         \right )
\end{align}
\label{thm:multinorms}
\end{theorem}
Moreover, there exists always an \emph{extremal multinorm} acting as storage function for internally stable and minimal systems, i.e. that attains the value $\gamma = \gamma_p(\graph, \mset)$ in (\ref{eq:optigain}).  
\begin{theorem}
Consider a system $(\graph(V,E), \mset)$, an integer $p \geq 1$, and $\gamma \geq \gamma_p(\graph, \mset)$. At each node $v \in V$, define the function

\begin{align}
F_{v,\gamma,p}(x) &= \s{\pi \in \graph, v_{\pi}(0) = v; \, \mathbf{w} \in \Ell_p }\|C_\pi x + D_\pi w\|_p^p - \gamma^p \|w\|_p^p, \nonumber\\
& = \sup_{\pi \in \graph, v_{\pi}(0) = v; \mathbf{w} \in \Ell_p} \sum \|z_t\|^p_p -  \gamma^p\sum \|w_t\|_p^p. \label{eq:Fv}
\end{align}
The set of functions $\{F_{v, \gamma, p}^{1/p}(x), \,  v \in V\}$ is an extremal multinorm, i.e., for all $(u,v,\sigma) \in E$, $x \in \reels^n$ and $w \in \reels^d$, 
\begin{equation}
\begin{aligned}
 F_{v,\gamma,p}(A_\sigma x + B_\sigma w) + \|C_\sigma x + D_\sigma w\|_p^p \\
 \leq F_{u,\gamma,p}(x) + (\gamma\|w\|_p)^p.
\end{aligned}
\label{eq:Fgamma}
\end{equation}
\label{thm:extremality}
\end{theorem}

\begin{shortVersion} 
\subsubsection*{Proofs of main results}
Consider a system $(\graph, \mset)$. We need to show that if $(\ref{eq:normineqgain})$ holds, then $\gamma \geq \gamma_p(\graph, \mset)$. Then, proving Theorem \ref{thm:extremality} proves Theorem \ref{thm:multinorms}. We know $\gamma_p(\graph, \mset)$ is bounded as a consequence of Assumption \ref{ass:stability}.\\
The first element is proven by simple algebra. Assume (\ref{eq:normineqgain}) holds, take $x_0 = 0$, an infinite path $\pi \in \graph$, and a sequence $\mathbf{w} \in \Ell_p$. Unfolding the inequality, we obtain for all $T \geq 1$,
$$
\begin{aligned}
 |x_0|_{v_{\pi}(0)}^p  & + \gamma^p (\|w_0\|_p^p + \sum_{t = 1}^\infty\|w_t\|_p^p) \\
 & \geq |x_1|_{v_{\pi}(1)}^p + \|z_0\|_{p}^p + \gamma^p(\sum_{t = 1}^\infty\|w_t\|_p^p)  \geq \cdots \\
 & \geq |x_T|_{v_{\pi}(T)}^p \hspace{-2pt} + \hspace{-2pt} \sum_{t = 0}^{T-1} \|z_t\|_p^p + \sum_{t = T}^\infty \|w_{t}\|_p^p \geq \sum_{t = 0}^{T-1} \|z_t\|_p^p.\\
 \end{aligned}$$
 Therefore, for all $T \geq 1$, 
 $ (\sum_{t = 0}^{T-1} \|z_t\|_p^p)/(\sum_{t = 0}^\infty \|w_t\|_p^p) \leq \gamma^p,$
 and this being independent of $\pi$ and $\mathbf{w} \in \mathcal{L}_p$, we get  $\gamma_p(\graph, \mset) \leq \gamma$.
We now move onto the proof of Theorem \ref{thm:extremality}. Proving (\ref{eq:Fgamma}) is easy, so we focus on proving that for $\gamma \geq \gamma_p(\graph, \mset)$,
the functions $F_{v,\gamma,p}^{1/p}$ are norms. We omit the subscripts $\gamma$ and $p$ in what follows. The functions have the following properties:
\begin{enumerate}
\item $\mathbf{F_{v}(0) = 0}.$  Observe that  $F_{v}(0) \geq 0$, taking a disturbance $\mathbf{w} = 0 \in \Ell_p$.
Assume by contradiction that there is $K > 0$ such that $F(0) \geq K$. There must be a path $\pi$ and a sequence $\mathbf{w} \in \mathcal{L}_p$ 
such that $$ \sum_{t = 0}^\infty \|z_t\|_p^p \geq K/2 + \gamma^p \sum_{t = 0}^\infty \|w_t\|_p^p > \gamma^p \sum_{t = 0}^\infty \|w_t\|_p^p. $$
Thus, $\gamma < \gamma_p(\graph, \mset)$ (see (\ref{eq:gain})), a contradiction.
\item $\mathbf{F_{v}(x) > 0}$ \textbf{is positive definite}. By Assumption \ref{ass:minimal}, for any node $v \in V$, for any $x \in \reels^n$, there is a path $\pi \in \graph$, $v_{\pi}(0) = v$, such that 
$C_\pi x \neq 0$. Thus,  $F_{v}(x) > 0$ is guaranteed by taking $\mathbf{w} = 0 \in \Ell_p$. 
\item $\mathbf{F_{v}(x)}$ \textbf{is convex and positively homogeneous of order} $\mathbf{p}$. 
Given a path $\pi \in \graph$ of finite length, define the output sequence
$$\mathbf{z}_\pi(x_0, \mathbf{w}) = \begin{pmatrix} z_0 \\ \vdots \\ z_{|\pi|-1} \end{pmatrix} = C_\pi x_0 + D_\pi \begin{pmatrix} w_0 \\ \vdots \\ w_{|\pi| - 1}\end{pmatrix}.$$
Clearly, $\mathbf{z}_\pi(\alpha x_0, \alpha \mathbf{w}) = \alpha \mathbf{z}_\pi(x_0, \mathbf{w})$ for any $\alpha \in \reels$. 
Therefore, for any $\alpha \in \reels$, $\alpha \neq 0$, we have
$$
\begin{aligned}
F_{v}( \alpha x) &= \s{\pi, \mathbf{w} \in \Ell_p} \sum_t \|z_t(\alpha x, \mathbf{w})\|_p^p - \gamma_p^p |w_t|_p^p, \\
& = \s{\pi, \alpha \mathbf{w'} \in \Ell_p} \sum_t  |\alpha|^p \|z_t(x, \mathbf{w'})\|_p^p - \gamma_p^p\alpha^p \|w'_t\|_p^p, 
\end{aligned} 
$$
Since $\|\mathbf{z}_\pi(x_0, \mathbf{w})\|_p^p = \sum_t \|z_t(x_0, \mathbf{w})\|_p^p$ is convex in $x_0$, the convexity of $F_v(x)$ is then easily proven. 
\item $\mathbf{F_v(x)< +\infty, \, \forall x \in \reels^n}$. 
Take  $x \in \im(B_\pi)$ for some $\pi \in \graph$ with $|\pi| < + \infty$, $v_{\pi}(|\pi|) = v$. Let $u = v_{\pi}(0)$, and let $\mathbf{w}$ be such that $ x = B_\pi \begin{pmatrix} w_0^\top, \ldots, w_{|\pi|-1}^\top \end{pmatrix} ^\top$. We have that 
$$
\begin{aligned}
F_{u}(0) = 0  \geq \sum_{t = 0}^{|\pi|-1}( \|z_t\|_p^p - \gamma_p^p \|w_t\|_p^p) + F_{v}(x), 
\end{aligned}
$$
and thus $F_{v}(x)$ is bounded.
Now, take any $x \in \reels^n$.
By Assumption $\ref{ass:minimal}$, we know that $x$ is a linear combination of points in the reachable sets $\im(B_\pi)$, $v_\pi(|\pi|) = v.$ Thus, by convexity, $F_v(x) < + \infty$. 
\end{enumerate}
\end{shortVersion}
We  discuss the generality of the Assumptions \ref{ass:minimal} and \ref{ass:stability} in the next subsections.
\subsection{Internal stability or undecidability}
\label{subsec:Stability}
Assumption \ref{ass:stability} is a sufficient condition for $\gamma_{p}(\graph, \mset) < + \infty$, which is key in Theorems \ref{thm:multinorms} and \ref{thm:extremality}. Relaxing the assumption leads to decidability issues:
\begin{proposition}
Given a switching system $(\graph, \mset)$ and $p \geq 1$, the question of whether or not $\gamma_p(\graph, \mset) < + \infty$ is \emph{undecidable}. 
\end{proposition}
\begin{shortVersion}
\begin{proof}
Consider an arbitrary switching system on two modes, with $\mathbf{A} = \{A_1, A_2\}$ and $\|A_i\| > 1$ for $i \in \{1,2\}$. The undecidability result presented in \cite{BlTsTBOA} states that the question of the existence of a uniform bound $K > 0$ such that, for all $T$ and all $\{\sigma(0), \ldots, \sigma(T-1)\} \in \{1,2\}^T$, the products of matrices satisfy $\|A_{\sigma(T-1)} ,\ldots, A_\sigma(0)\| \leq K$, is undecidable.
We show that as a consequence, one cannot in general decide the boundedness of the gain of a (constrained) switching system. Given  a pair of matrices, construct an arbitrary switching system of the form (\ref{eq:switching}) on 3 modes with the parameters
$$\mset = \{\{A_1,0,0,0\}, \{A_2,0,0,0\}, \{0,\mathbf{I},\mathbf{I},0\}\}, $$
with input and output dimension $d =  m = n$.  
Then, the $\mathcal{L}_p$-gain of this system is given by 
$$\gamma_p(\graph, \mset) = \hspace{-4pt} \s{\substack{w_0, T, \\ \sigma(0), \ldots, \sigma(T-1) \in \{1,2\}^T}} \hspace{-4pt} \frac{\|A_{\sigma(T-1)} \cdots A_{\sigma(0)}w_0\|_p}{\|w_0\|_p}. $$
Thus, the gain of the system  $(\graph, \mset)$ is bounded if and only if there is a uniform bound on the norm of all products of matrices from the pair $\mathbf{A} = \{A_1, A_2\}$. 
\end{proof}
\end{shortVersion}
\subsection{Obtaining minimal realizations.}
\label{subsec:Minimize}
We now turn our attention to Assumption \ref{ass:minimal}. \\
Section \ref{sec:Example} presents a practically motivated system for which a very natural model turns out to be non-minimal. We need to provide a minimal realization for the system. Algorithms exist in the literature for computing minimal realizations for arbitrary switching systems (see \cite{PeBaRTOD}). We use an approach similar to that of \cite{PeBaRTOD}, for constrained switching systems. This produces a system with the same input-output relations, but that is \emph{rectangular}, with a dimension of the state space varying in time.   
\begin{definition} 
A \emph{rectangular} constrained switching system a tuple $(\graph(V,E), \mset, \mathbf{n})$.  The graph $\graph$ is strongly connected with nodes $V$ and labeled, directed edges $E$. Edges are of the form $(v_1, v_2, \sigma) \in E$, where $\sigma \in \{1, \ldots, |E|\}$ are \emph{unique} labels assigned to each edge. The \emph{dimension vector} $\mathbf{n} = (n_v)_{v \in V}$ assigns a \emph{dimension} to each node $v \in V$. The set $\mset = \{\{A_\sigma, B_\sigma, C_\sigma, D_\sigma\}_{\sigma \in 1, \ldots, |E|}\}$ satisfies, for $(u, v, \sigma) \in E$, 
$A_\sigma \in \reels^{n_{v} \times n_{u}}, \, B_{\sigma} \in \reels^{n_{v} \times d}, \, C_{\sigma} \in \reels^{m \times n_{u}}, \, D_\sigma \in \reels^{m \times d}. $
$\,$
\end{definition}
\begin{remark}
The systems considered so far can be cast as rectangular systems, with $\mathbf{n} = \{n_v = n, v \in V\}$. The one-to-one correspondence between modes and edges is easily obtained through relabeling.
\label{rem:torect}
\end{remark}
\begin{remark}
Given any pair  $(v_1, v_2, \sigma_1), (v_2,v_3,\sigma_2) \in E$, the products $A_{\sigma_2} A_{\sigma_1}$, $C_{\sigma_2}A_{\sigma_1}$, etc... are compatible.
Thus, for a rectangular system $(\graph, \mset, \mathbf{n})$ and a path $\pi \in \graph$, we may still compute the matrices 
$A_\pi, \, B_\pi, \, C_\pi$ and $D_\pi$ (see (\ref{eq:Amap}), (\ref{eq:Bmap}), (\ref{eq:Cmap}), (\ref{eq:Dmap})). The definitions for the subspaces $\mathcal{B}_v \subseteq \reels^{n_v}$ and $\mathcal{C}_v \subseteq \reels^{n_v} $ (see Assumption \ref{ass:minimal}) for rectangular systems  follow immediately.
\end{remark}
\begin{definition}
A rectangular system $(\graph(V,E), \mset, \mathbf{n})$, $\mathbf{n}= (n_v)_{v \in V}$, is said to be \emph{minimal} if, for all $v \in V$, 
$\mathcal{B}_v = \reels^{n_v},$ and $\mathcal{C}_v = \{0\} \subset \reels^{n_v}. $
\end{definition}
\begin{remark}
Rectangular switching systems are only special cases of more general hybrid systems. The results we present here regarding the minimization procedure, and the mechanisms behind them, can be deduced from those presented in \cite{PeScRTOD,PeBaRTOD}. 
\end{remark}
In order to compute a minimal realization for a system $(\graph, \mset, \mathbf{n})$ we first compute the subspaces $\mathcal{B}_v$ and $\mathcal{C}_v$. 
\begin{proposition}[Construction of the subspaces $\mathcal{C}_v$]
Given a rectangular system $(\graph, \mset, \mathbf{n})$, let for $v \in V$
$$ X_{v,1} = \sum_{(v,u,\sigma) \in E} C_\sigma^\top C_\sigma; \, \mathcal{C}_{v,1} = \ker(X_{v,1}). $$
For $k = 1,2, \ldots$, consider the following iteration
$$ X_{v,k+1} = \sum_{(v,u,\sigma) \in E} A_\sigma^\top X_{u,k} A_\sigma; \, \mathcal{C}_{v,k} = \ker( \sum_{t = 1}^k X_{v,k}). $$
The sequence $\mathcal{C}_{v,k}$ converges and $\mathcal{C}_{v} = \mathcal{C}_{v,K}$, $K = \sum n_v$. 
\label{prop:algoForC}
\end{proposition}
\begin{shortVersion}
\begin{proof}
We construct $X_{v,k}$ in such a way that for $x \in \reels^n$, $X_{v,l}x = 0$ if and only if $x \in \mathcal{C}_{v,k}$, that is, for any path initiated at $v$ and of length $k$, $x$ generates only the 0 output. 
The argument for convergence is that if at some iteration, $\mathcal{C}_{v,k+1} = \mathcal{C}_{v,k}$ for all $v \in V$, then the algorithm terminates. Thus, there can at most $K = \sum n_v$ steps.
\end{proof}
\end{shortVersion}
\begin{remark}
We may use the procedure above to construct $\mathcal{B}_v$. Indeed,
we can write $$\mathcal{B}_v^\perp = \bigcap_{\pi \in \graph; \, v_\pi(|\pi|) = v} \ker(B_\pi^\top). $$
Thus, given $(\graph(V,E), \mset, \mathbf{n})$ it suffice to apply Proposition \ref{prop:algoForC} to the \emph{dual system} $(\graph(V,\bar{E}), \bar{\mset}, \mathbf{n})$, defined with
$\bar{E} = \{(v,u,\sigma) : (u,v,\sigma) \in E \}$, and 
$$\bar{\mset} = \{(A_\sigma^\top, C_\sigma^\top, B_\sigma^\top, D_\sigma^\top): \,(A_\sigma, B_\sigma, C_\sigma, D_\sigma) \in \mset \},$$ to compute $\mathcal{B}_v^\perp $.  
\end{remark}

\begin{proposition}
Given a rectangular system $(\graph, \mset,\mathbf{n})$, let $n_v$ be the dimension of the node $v \in V$ and $m_v \leq n_v$ be the dimension of $\mathcal{B}_v$. 
For $v \in V$, let $L_v \in \reels^{n_v \times m_v}$ be an orthogonal basis of $\mathcal{B}_v$. The rectangular system $(\graph(V,E), \bar{\mset}, \bar{\mathbf{n}})$ on the set 
$$\bar{\mset} = \left \{
 (
\begin{aligned}
L_v^\top A_\sigma L_u , \,
L_v^\top B_\sigma, \,
 C_\sigma L_u , \,
D_\sigma    ,
    \end{aligned}
    )_{(u,v,\sigma) \in E} \right \} ,$$
has $\mathcal{B}_v = \reels^{m_v}$, $\bar{\mathbf{n}} = \{m_v, v\in V\}$, and $\gamma_p(\graph,\mset, \mathbf{n}) = \gamma_p(\graph, \bar{\mset}, \bar{\mathbf{n}})$.
\label{prop:getBvright}
\end{proposition}
\begin{shortVersion}
\begin{proof}
Observe that the subspaces $\mathcal{B}_v$ are invariant in the sense that for any $(u,v,\sigma) \in E$, 
$ \mathcal{B}_v \supseteq A_\sigma \mathcal{B}_u + \im(B_\sigma).$
Thus, if $x_0 = 0$, we know that after $t$ steps, if the next edge is $(u,v,\sigma)$, 
$$
\begin{aligned}
L_v^\top L_vx_{t+1} &= A_\sigma L_u^\top L_u x_t + B_\sigma w_t,\\
z_t &= C_\sigma L_u^\top L_u x_t + D_\sigma w_t.
\end{aligned}
$$
The system on the parameters $\bar{\Sigma}$ is then obtained by multiplying the first equation on the left by $L_v$ (since $L_vL_v^\top L_v = L_v$), and then doing the change of variable $y_t = L_{v_{\pi}(t)}^\top x_t$.
\end{proof}
\end{shortVersion}
\begin{proposition}
Given a rectangular system $(\graph, \mset,\mathbf{n})$, let $n_v$ be the dimension of the node $v \in V$ and $m_v \leq n_v$ be the dimension of $\mathcal{C}_v^\perp$. 
Let $K_v \in \reels^{n_v \times m_v}$ be an orthogonal basis of $\mathcal{C}_v^\perp$. The rectangular system $(\graph(V,E), \bar{\mset}, \bar{\mathbf{n}})$ on the set 
$$\bar{\mset} = \left \{
 (
\begin{aligned}
K_v^\top A_\sigma K_u , \,
K_v^\top B_\sigma, \,
 C_\sigma K_u , \,
D_\sigma    ,
    \end{aligned}
   )_{(u,v,\sigma) \in E}  \right \} ,$$
has $\mathcal{C}_v^\perp = \reels^{m_v}$, $\bar{\mathbf{n}} = \{m_v, v \in V\}$, and $\gamma_p(\graph,\mset, \mathbf{n}) = \gamma_p(\graph, \bar{\mset}, \bar{\mathbf{n}})$.
\label{prop:getCvright}
\end{proposition}
\begin{shortVersion}
\begin{proof}
The proof is similar to that of Proposition \ref{prop:getCvright}. Observe that  for any $(u,v,\sigma) \in E$, 
$ \mathcal{C}_v \supseteq A_\sigma \mathcal{C}_u.$
What we do then is to project the state onto $\mathcal{C}_v^\perp$ at every step. Doing so does not affect the gain since the part of the state in $\mathcal{C}$ would not affect the output at any time after the projection (by invariance).
\end{proof}
\end{shortVersion}
\begin{theorem}
Given a rectangular system $(\graph, \mset, \mathbf{n}= \{n_v, v \in V\})$,   consider the following iteration procedure. Let $S_0 = (\graph, \mset, \mathbf{n})$, and for $k = 1,2,\ldots,$ let $S_k$ be the system obtained by applying first Proposition \ref{prop:getBvright} then  Proposition \ref{prop:getCvright} to $S_{k-1}$. 
Then  $S_{N}$   is minimal for $N = \sum n_v$, and $\gamma_p(\graph, \mset, \mathbf{n}) = \gamma_p(\graph, \bar{\mset}, \bar{\mathbf{n}})$. 
\label{thm:minimize}
\end{theorem}
\begin{shortVersion}
\begin{proof}
The fact that the gain is conserved along the iterations is granted from Proposition \ref{prop:getBvright} and Proposition \ref{prop:getCvright}. The number of iterations needed to converge is actually conservative, but is easily shown. If the system is non-minimal, then at least one node will have a reduced dimension after an iteration. We can register a decrease at at least one node at most   $\sum_{v} n_v$ times.
\end{proof}
\end{shortVersion}

As a conclusion, even when given a non-minimal system $(\graph, \mset)$,   we may always construct a minimal \emph{rectangular} system  $(\graph, \bar{\mset}, \{(n_v)_{v\in V}\})$ with the same $\Ell_p$-gain.
Theorem \ref{thm:multinorms} is easily translated for rectangular systems, by defining a multinorm $\{(|{\cdot}|_v)_{v \in V}\}$ 	as a set of norms, where $|{\cdot}|_v$ is a norm on the space $\reels^{n_v}$.  The storage functions of Theorem \ref{thm:extremality} remain well defined, with  $F_{v, \gamma,p}$ now taking values  in $\reels^{n_v}$. 
 \begin{remark}
There might be some pathological cases where the dimension of some nodes reduces to $0$. This would correspond to a situation where, initially, $\mathcal{B}_v \subseteq \mathcal{C}_v$.  
\end{remark}

\section{Storage functions for the $\Ell_2$-gain.}
 \label{sec:L2}
We now aim at computing the extremal storage functions $F_{v, \gamma,2}(x)$ (see Theorem \ref{thm:extremality}) of a system for computing its $\Ell_2$-gain. We no longer consider rectangular systems here, but the generalization of the results is straightforward. The Assumptions \ref{ass:minimal} and \ref{ass:stability} still hold. 
 
In approximating the function $F_{v, \gamma,2}(x)$, two questions arise. First, given a value of $\gamma$, how can we get a good approximation of $F_{v, \gamma,2}(x)$? Second, and  more importantly, how can we obtain obtain a good estimate of the $\Ell_2$-gain? \\
We propose a first method, based on dynamic programming and asymptotically tight under-approximations of the  $\Ell_2$-gain, and a second method, based on the path-dependent Lyapunov function framework of \cite{EsLeCOLS,LeDuUSOD}, for obtaining over-approximations of the $\Ell_2$-gain. 

Both are based on the following observation. Given a system $(\graph(V,E), \mset)$, we can write  
$$
F_{v, \gamma,p}(x) = \s{\pi \in \graph, |\pi| = \infty, v_{\pi}(0) = v} F_{v, \gamma,p, \pi}(x),
$$
with
\begin{equation}
F_{v, \gamma,p, \pi}(x) = \s{\mathbf{w} \in \Ell_p }\|C_\pi x + D_\pi w\|_p^p - \gamma^p \|w\|_p^p, \label{eq:Fgpi}
\end{equation}
where the path $\pi$ is fixed. The definition above holds for any path of  finite or infinite length. A first approximation of $F_{v, \gamma, p}$ is obtained limiting the length of the paths to some $K \geq 1$:
$$\check{F}_{v, \gamma, p,K} = \max_{\pi \in \graph, v_{\pi}(0)= v, |\pi| = K}F_{v, \gamma, p, \pi}.$$ 
We have the following for $p = 2$:
\begin{align}
& F_{v, \gamma,2, \pi}  (x) = \label{eq:F2} \\
&   x^\top \left( C_\pi^\top C_\pi - C_\pi^\top D_\pi (D_\pi^\top D_\pi - \gamma^2 {I}_d) ^{-1}  D_\pi^\top C_\pi \right) x. \nonumber
\end{align}
If $|\pi| = K$, $F_{v, \gamma, 2, \pi}$ can be computed using dynamic programming. Indeed,
$$\begin{aligned}
F_{v, \gamma, 2, \pi}(x)& \\
& \hspace{-28pt} = \max_{\mathbf{w}} \|C_{\pi(1:K-1)}x \hspace{-1pt} + \hspace{-1pt} D_{\pi(1:K-1)} (w_0^\top, \ldots, w_{K-2}^\top)^\top \hspace{-1pt}\|^2_2 \\
& \hspace{-18pt} + \max_{w_{K-1}} \|C_{\sigma_{\pi}(K)} A_{\pi(1:K-1)} x + D_{\sigma_{\pi}(K)} w_{K-1}\|^2_2,
\end{aligned} $$
 We can then solve for $w_{K-1}$ as a function of $y = A_{\pi(1:K-1)}x$, and then solve for $w_{K-2}, \ldots, w_0$ in succession.
\begin{proposition}
Given a system $(\graph(V,E), \mset)$, and an integer $K \geq 1$ let 
$$\check{\gamma}_{K,p} = \s{\pi \in \graph, |\pi| = K} \|D_\pi\|_p.$$ 
For any $\gamma > \check{\gamma}_{K,2}$, $K \geq n |V|$, $\check{F}_{v, \gamma , 2,K}^{1/2}$ is a norm.
\label{prop:truncatedNorms}
\end{proposition}
\begin{shortVersion}
\begin{proof} 
The term $(D_\pi^\top D_\pi - \gamma^2 I)$ in (\ref{eq:F2}) is by definition negative definite if $\gamma > \check{\gamma}_{K,2}$.
Thus, $F_{v, \gamma, 2, \pi}$ is well defined. Then, by Assumption \ref{ass:minimal} and Proposition \ref{prop:algoForC}, for  $K \geq n|V|$, for all $x \in \reels^n$, there is $\pi$ such that $C_\pi x \neq 0$.  We can then conclude that $\check{F}^{1/2}$ is a norm. 
\end{proof}
\end{shortVersion}
\begin{proposition}
Given a system $(\graph(V,E), \mset)$, for any $k \geq 1$, $\check{\gamma}_{k,p} \leq \gamma_p(\graph, \mset)$, and $\lim_{k \rightarrow \infty} \check{\gamma}_{k,p} =\gamma_p(\graph, \mset)$.
\label{prop:gainLowerbound}
\end{proposition}
The approach above allows to get asymptotically tight lower-bounds. 
In order to get upper bounds, we can \emph{approximate} storage functions by quadratic norms using  the Horizon-Dependent Lyapunov functions of \cite{EsLeCOLS}.
\begin{proposition}[\cite{EsLeCOLS}]
Consider a system $(\graph, \mset)$. For any $K \geq 1$, consider the following program:
$$ 
\begin{aligned}
& \hat\gamma_K= \inf_{\gamma,\, X_\pi \in \reels^{n \times n}, \, \pi \in \graph, \, |\pi| = K } \gamma \qquad s.t.\\
& \forall \pi_1, \pi_2 \in \graph, \, |\pi_1| \hspace{-1pt} = \hspace{-1pt} |\pi_2| \hspace{-1pt} = \hspace{-1pt} K, \, \pi_1(2:K) = \pi_2(1:K-1), \\
& 
\begin{pmatrix}
A_{\sigma_{\pi_1}(0)} & B_{\sigma_{\pi_1}(0)} \\
D_{\sigma_{\pi_1}(0)}  & C_{\sigma_{\pi_1}(0)} 
\end{pmatrix}^\top 
\begin{pmatrix}
X_{\pi_2} & 0 \\
0 & I
\end{pmatrix}
\begin{pmatrix}
A_{\sigma_{\pi_1}(0)} & B_{\sigma_{\pi_1}(0)} \\
D_{\sigma_{\pi_1}(0)}  & C_{\sigma_{\pi_1}(0)} 
\end{pmatrix}\\
& \qquad  -
\begin{pmatrix}
X_{\pi_1} & 0 \\
0 & \gamma^2I
\end{pmatrix} \preceq 0,\\
&\forall \pi \in \graph, \, |\pi| = K, \, X_\pi \succ 0. 
\end{aligned}
$$
Then $\lim_{K \rightarrow \infty} \hat\gamma_K = \gamma_p(\graph, \mset),$ and $\hat \gamma_K \geq \gamma_p(\graph, \mset).$
\label{prop:gainUpperbound}
\end{proposition}
\begin{remark}
There is an interesting link between Horizon-Dependent Lyapunov functions and the approximation of the functions $F_{v,\gamma,p}$.
Consider the function 
 $$
\hat{F}_{v, \gamma,p, \pi}(x) = \s{\phi \in \graph, \, \phi(1:K) = \pi(1:K)} F_{v, \gamma, p, \phi}(x),
$$
with $F_{v, \gamma, p, \phi}$ as in (\ref{eq:Fgpi}). Compared to $F_{v, \gamma, p}$, we now fix the first $K$ edges of the path considered in (\ref{eq:Fv}) to be those of a path $\pi$. It is easily seen that
$$
F_{v, \gamma,p}(x) = \max_{\pi \in \graph, v_{\pi}(0) = v, |\pi| = K} \hat{F}_{v, \gamma,p, \pi}(x).
 $$
 By definition, if we take two paths $\pi_1$ and $\pi_2$ of length $K$ such that $\pi_1(2:K) = \pi_2(1:K-1)$, we can then verify that
$$ \begin{aligned}
\hat{F}_{v, \gamma,p, \pi_1}(x) + \gamma^p \|w\|^2_2 & \geq \hat{F}_{v, \gamma,p, \pi_2}(A_{\sigma_{\pi}(0)}x + B_{\sigma_{\pi}(0)}w) \\
&+  \|C_{\sigma_{\pi}(0)}x + D_{\sigma_{\pi}(0)}w)\|^2_2.
\end{aligned}
$$
For $p = 2$, assuming the functions $\hat{F}$ to be quadratic, the inequalities above are equivalent to the LMIs of Proposition \ref{prop:gainUpperbound}.
\end{remark}

Together, Propositions \ref{prop:gainLowerbound} and \ref{prop:gainUpperbound} enable us to approximate the $\Ell_2$-gain of a system, and its storage functions, in an arbitrarily accurate manner.

\subsection{Converse results for the existence of quadratic storage functions.}
\label{Subsection:Converse}

Quadratic (multiple) Lyapunov functions have received a lot of attention in the past for the stability analysis of switching systems (see e.g. \cite{LiMoBPIS,JuTJSR, BlFeSAOD, PhEsSODT}). Checking for their existence is computationally easy as it boils down to solving LMIs. They are however conservative certificates of stability, and thus it is interesting to seek ways to quantify \emph{how conservative} these methods are. In the following, we extend existing results on the conservatism of quadratic Lyapunov functions (see \cite{AnShSC}, \cite{BlNeOTAO}, \cite{PhEsSODT}) to the performance analysis case. We give a converse theorem for the existence of quadratic storage functions, along with a conjecture related to Horizon-Dependent Storage functions (Proposition \ref{prop:gainUpperbound}).
\begin{theorem}
Given a constrained switching system $(\graph(V,E), \mset)$, consider the system $(\graph, \mset')$ with
$$\mset' = \{ (\sqrt{n}A_\sigma, \, \sqrt{n}B_\sigma, \, C_\sigma, D_\sigma): \,(A_\sigma, \, B_\sigma, \, C_\sigma, D_\sigma) \in \mset\}.$$
If $\gamma_{2}(\graph, \mset') \leq 1$, then $\gamma_2(\graph, \mset) \leq 1$ and there is a set of \emph{quadratic} norms $\{(|\cdot|_{Q,v})_{v\in V}\}$ such that for all 
$x \in \reels^n$, $w \in \reels^d$ and $(u,v,\sigma) \in E$,
$$|A_\sigma x + B_\sigma w|_{Q,v}^{2} + \|C_\sigma x + D_\sigma w\|_2^2 \leq |x|_{Q,u}^{2} +  \|w\|_2^{2}.$$
\end{theorem}
\begin{shortVersion}
\begin{proof}
The result relies on Theorem \ref{thm:multinorms}, more precisely on the fact that the functions $F_{v, \gamma, 2}$ are norms. The proof is similar to that of \cite{PhEsSODT}, Theorem 3.1: we can use John's ellipsoid theorem (see e.g. \cite{BlNeOTAO}) to approximate the norms of a storage function for $(\graph,\mset')$ with quadratic norms. These quadratic norms then provide a storage function for $(\graph, \mset)$.
\end{proof}
\end{shortVersion}
We conjecture that the following extension, which follows from the stability analysis case (see e.g. \cite{PhEsSODT}, Theorem 3.5), holds true  for the performance analysis case:
\begin{conjecture}
Given a constrained switching system $(\graph(V,E), \mset)$, consider the system $(\graph, \mset')$ with
$$\mset' = \{ (n^{\frac{1}{2d}}A_\sigma, \, n^{\frac{1}{2d}}B_\sigma, \, C_\sigma, D_\sigma): \,(A_\sigma, \, B_\sigma, \, C_\sigma, D_\sigma) \in \mset\}.$$
If $\gamma_{2}(\graph, \mset') \leq 1$, then $(\graph, \mset)$ has a Horizon-Dependent storage function (see Proposition \ref{prop:gainUpperbound}) for $K = d+1$.
\end{conjecture}

\section{Example.}
\label{sec:Example}
\begin{figure}[!ht]
\centering
 \includegraphics[width = 0.4\columnwidth]{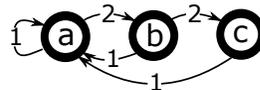}
 \caption{Graph for the example. Here, we cannot have more than 2 failures (that is, mode $\sigma(t) = 2$) in a row. Nodes have been labeled $a$,$b$, and $c$ for further discussions.}
 \label{fig:graphExample}
 \end{figure}

We are given a stabilized LTI system\footnote{We use the simple model of an inverted pendulum with mass $2$kg. The system is linearized around the ``up'' position, and discretized at 100 hz. The control gains are computed through LQR, with cost 1 on the norm of the output, and 10 on the norm of the input. Computations done in Matlab, codes available at \url{http://sites.uclouvain.be/scsse/gainsAndStorage.zip}}, 
$$x_{t+1} = A x_t + B u_t, \, z_t =  x_t, \, u_t = K x_t. $$
We let $x \in \reels^2$, $z \in \reels^2$, $w \in \reels^1$.
We assume that there might be \emph{delays} in the control updates. At any time, either $u_t = Kx_t$, or $u_t = u_{t-1}$, and we assume that there cannot be more than two delays in a row. Moreover, when there is a delay,  the system undergoes disturbances from the actuator, i.e. of the form $Bw_t$. 
The situation can be modeled as a constrained switching system  $(\graph, \mset)$ on two modes (with $\sigma(t) = 1$ if the control input is updated or $\sigma(t) = 2$ else) with the graph $\graph$ of Figure \ref{fig:graphExample} and
$$\mset =  \left \{ 
\begin{aligned}
&\left (  \begin{pmatrix} A+BK & 0 \\ I & 0  \end{pmatrix}, \, \begin{pmatrix} 0 \\  0   \end{pmatrix}, \, 
\begin{pmatrix} I &  0  \end{pmatrix},  \, \begin{pmatrix} 0 \end{pmatrix} \right )\\  
&\left ( \begin{pmatrix} A & BK  \\ 0 & I  \end{pmatrix}, \, \begin{pmatrix} B \\  0   \end{pmatrix}, \, 
\begin{pmatrix} I &  0   \end{pmatrix},  \, \begin{pmatrix} 0 \end{pmatrix} \right )
\end{aligned}
 \right   \}.$$
 \vspace{-10pt}
 
\begin{figure}[!ht]
\centering
\begin{subfigure}[c]{\columnwidth}
\centering
\includegraphics[scale = 0.35]{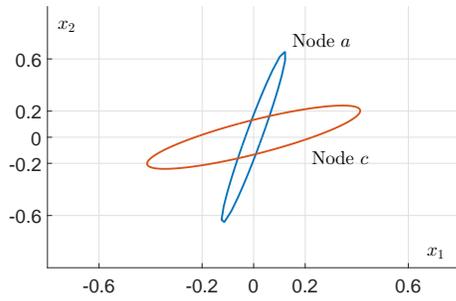}
\caption{Level sets of $\check{F}_{a,\tilde{\gamma},2}^{1/2}(x)$ and $\check{F}_{c,\tilde{\gamma},2}^{1/2}(x)$.}
\label{fig:storageA}
\end{subfigure}
\begin{subfigure}[c]{\columnwidth}
\centering
\includegraphics[scale = 0.35]{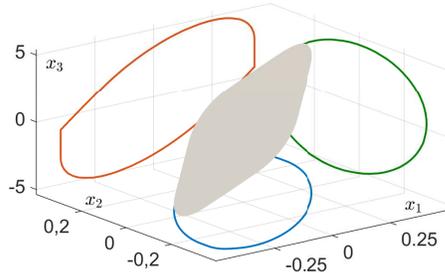}
\caption{Level sets of $\check{F}_{b,\tilde{\gamma},2}(x)$, with projection on each coordinate plane.}
\label{fig:storageB}
\end{subfigure}
\caption{Level sets of the storage functions approximation at each node. Node $b$ has a state dimension of 3.}
\label{fig:norms}
\end{figure}

The system here is actually \emph{not minimal}. Indeed, take the third node $c$. Any vector of the form $ x =\begin{pmatrix}0 & 0 & \alpha & \beta \end{pmatrix} ^\top$ satisfies $C_\pi x = 0 $, for any $\pi$ such that $v_\pi(0) = c$. Thus, $\dim(\mathcal{C}_c) = 2$. From there, we can see that at nodes $a$ and $b$, $\mathcal{C}_a$ and $\mathcal{C}_b$ are formed of any vector of the form $ x =\begin{pmatrix}0 & 0 & \alpha & \beta \end{pmatrix}^\top$, where $BK\begin{pmatrix} \alpha & \beta \end{pmatrix}^\top = 0$. Thus, $\dim(\mathcal{C}_a) = \dim(\mathcal{C}_b) =  1$. Applying Theorem \ref{thm:minimize} to the system $(\graph, \mset)$ (after applying Remark \ref{rem:torect}), we obtain a minimal rectangular system with nodal dimensions $(2,3,2)$. We can then approximate the storage functions $F_{v,\gamma,2}^{1/2}, v \in \{a,b,c\}$.
 Applying Proposition \ref{prop:truncatedNorms} with paths of length 10, we obtain, a  lower-bound on the $\Ell_2$ gain
$$ \tilde{\gamma}= \sup_{\pi \in \graph, |\pi| = K} \|D_\pi\|_p  \simeq 0.0188. $$
We use this bound to compute the approximations of $F_{v,\gamma,2}^{1/2}, v \in \{a,b,c\}$, $\check{F}_{v,\tilde{\gamma},2}^{1/2}, v \in \{a,b,c\}$. The level sets of these functions are displayed on Figure \ref{fig:norms}. 

\section{Conclusion.}

We provide a general characterization of the $\Ell_p$-gain of discrete-time linear switching system under the form of switching storage functions. Under the assumptions of internal stability and minimality, the pieces of these functions are the $p$th power of norms.  The generality of these assumptions is discussed and we provide means to compute minimal realizations of constrained switching systems. We then turn our focus on the $\Ell_2$ gain, and provide algorithms for obtaining  asymptotically tight lower and upper bounds on the gain based on the approximation of storage functions. Finally, we provide a converse result for the existence of quadratic storage functions exploiting the nature of storage functions, and formulate a conjecture about Horizon-Dependent storage functions. We believe an answer to the conjecture (positive or negative) will allow for a better understanding of the geometry underlying Lyapunov methods for performance analysis.

\bibliographystyle{ieeetrans}

\bibliography{biblio}

\begin{thebibliography}{10}
\providecommand{\url}[1]{#1}
\csname url@samestyle\endcsname
\providecommand{\newblock}{\relax}
\providecommand{\bibinfo}[2]{#2}
\providecommand{\BIBentrySTDinterwordspacing}{\spaceskip=0pt\relax}
\providecommand{\BIBentryALTinterwordstretchfactor}{4}
\providecommand{\BIBentryALTinterwordspacing}{\spaceskip=\fontdimen2\font plus
\BIBentryALTinterwordstretchfactor\fontdimen3\font minus
  \fontdimen4\font\relax}
\providecommand{\BIBforeignlanguage}[2]{{%
\expandafter\ifx\csname l@#1\endcsname\relax
\typeout{** WARNING: IEEEtranS.bst: No hyphenation pattern has been}%
\typeout{** loaded for the language `#1'. Using the pattern for}%
\typeout{** the default language instead.}%
\else
\language=\csname l@#1\endcsname
\fi
#2}}
\providecommand{\BIBdecl}{\relax}
\BIBdecl

\bibitem{AhJuJSRA}
A.~A. Ahmadi, R.~M. Jungers, P.~A. Parrilo, and M.~Roozbehani, ``Joint spectral
  radius and path-complete graph lyapunov functions,'' \emph{SIAM Journal on
  Control and Optimization}, vol.~52, no.~1, pp. 687--717, 2014.

\bibitem{AlRaCMAA}
R.~Alur, A.~D'Innocenzo, K.~H. Johansson, G.~J. Pappas, and G.~Weiss,
  ``Compositional modeling and analysis of multi-hop control networks,''
  \emph{IEEE Transactions on Automatic control}, vol.~56, no.~10, pp.
  2345--2357, 2011.

\bibitem{AnShSC}
T.~Ando and M.-H. Shih, ``Simultaneous contractibility,'' \emph{SIAM Journal on
  Matrix Analysis and Applications, 19(2), 487-498}, 1998.

\bibitem{BlFeSAOD}
P.-A. Bliman and G.~Ferrari-Trecate, ``Stability analysis of discrete-time
  switched systems through lyapunov functions with nonminimal state,'' in
  \emph{Proceedings of IFAC Conference on the Analysis and Design of Hybrid
  Systems}, 2003, pp. 325--330.

\bibitem{BlNeOTAO}
V.~D. Blondel, Y.~Nesterov, and J.~Theys, ``On the accuracy of the ellipsoid
  norm approximation of the joint spectral radius,'' \emph{Linear Algebra and
  its Applications}, vol. 394, pp. 91--107, 2005.

\bibitem{BlTsTBOA}
V.~D. Blondel and J.~N. Tsitsiklis, ``The boundedness of all products of a pair
  of matrices is undecidable,'' \emph{Systems \& Control Letters}, vol.~41,
  no.~2, pp. 135--140, 2000.

\bibitem{CoBoRMSG}
P.~Colaneri, P.~Bolzern, and J.~C. Geromel, ``Root mean square gain of
  discrete-time switched linear systems under dwell time constraints,''
  \emph{Automatica}, vol.~47, no.~8, pp. 1677--1684, 2011.

\bibitem{DaBePDLF}
J.~Daafouz and J.~Bernussou, ``Parameter dependent lyapunov functions for
  discrete time systems with time varying parametric uncertainties,''
  \emph{Systems \& Control Letters}, vol.~43, no.~5, pp. 355--359, 2001.

\bibitem{DaAGTS}
X.~Dai, ``A {G}el'fand-type spectral radius formula and stability of linear
  constrained switching systems,'' \emph{Linear Algebra and its Applications},
  vol. 436, no.~5, pp. 1099--1113, 2012.

\bibitem{DuLaANAF}
G.~E. Dullerud and S.~Lall, ``A new approach for analysis and synthesis of
  time-varying systems,'' \emph{IEEE Transactions on Automatic Control},
  vol.~44, no.~8, pp. 1486--1497, 1999.

\bibitem{EsLeCOLS}
R.~Essick, J.-W. Lee, and G.~E. Dullerud, ``Control of linear switched systems
  with receding horizon modal information,'' \emph{IEEE Transactions on
  Automatic Control}, vol.~59, no.~9, pp. 2340--2352, 2014.

\bibitem{EsPhTMAS}
R.~Essick, M.~Philippe, G.~Dullerud, and R.~M. Jungers, ``The minimum
  achievable stability radius of switched linear systems with feedback,'' in
  \emph{IEEE 54th Annual Conference on Decision and Control (CDC)}.\hskip 1em
  plus 0.5em minus 0.4em\relax IEEE, 2015, pp. 4240--4245.

\bibitem{HeMiOAMS}
E.~A. Hernandez-Vargas, R.~H. Middleton, and P.~Colaneri, ``Optimal and mpc
  switching strategies for mitigating viral mutation and escape,'' in
  \emph{Proc. of the 18th IFAC World Congress Milano (Italy) August}, 2011.

\bibitem{HeLIGO}
J.~P. Hespanha, ``L2-induced gains of switched linear systems,'' \emph{Unsolved
  problems in mathematical systems and control theory}, p. 131, 2004.

\bibitem{HiHeLIGA}
K.~Hirata and J.~P. Hespanha, ``L2-induced gain analysis of switched linear
  systems via finitely parametrized storage functions,'' in \emph{American
  Control Conference (ACC), 2010}.\hskip 1em plus 0.5em minus 0.4em\relax IEEE,
  2010, pp. 4064--4069.

\bibitem{JuTJSR}
R.~Jungers, ``The joint spectral radius,'' \emph{Lecture Notes in Control and
  Information Sciences}, vol. 385, 2009.

\bibitem{JuDiFSOD}
R.~M. Jungers, A.~D'Innocenzo, and M.~D. Di~Benedetto, ``Feedback stabilization
  of dynamical systems with switched delays,'' in \emph{Proc. of the 51st IEEE
  Conference on Decision and Control}, 2012, pp. 1325--1330.

\bibitem{KoTBWF}
V.~Kozyakin, ``The {B}erger--{W}ang formula for the markovian joint spectral
  radius,'' \emph{Linear Algebra and its Applications}, vol. 448, pp. 315--328,
  2014.

\bibitem{KuChSSSF}
A.~Kundu and D.~Chatterjee, ``Stabilizing switching signals for switched
  systems,'' \emph{IEEE Transactions on Automatic Control,}, vol.~60, no.~3,
  pp. 882--888, 2015.

\bibitem{LeDuODAF}
J.-W. Lee and G.~E. Dullerud, ``Optimal disturbance attenuation for
  discrete-time switched and markovian jump linear systems,'' \emph{SIAM
  Journal on Control and Optimization}, vol.~45, no.~4, pp. 1329--1358, 2006.

\bibitem{LeDuUSOD}
------, ``Uniform stabilization of discrete-time switched and markovian jump
  linear systems,'' \emph{Automatica, 42(2), 205-218}, 2006.

\bibitem{LiSISA}
D.~Liberzon, \emph{Switching in systems and control}.\hskip 1em plus 0.5em
  minus 0.4em\relax Springer Science \& Business Media, 2012.

\bibitem{LiMoBPIS}
D.~Liberzon and A.~S. Morse, ``Basic problems in stability and design of
  switched systems,'' \emph{IEEE Control Systems Magazine}, vol.~19, no.~5, pp.
  59--70, 1999.

\bibitem{NaVoCAOO}
M.~Naghnaeian and P.~Voulgaris, ``Characterization and optimization of l1 gains
  of linear switched systems.''

\bibitem{NaVoAUFF}
M.~Naghnaeian, P.~G. Voulgaris, and G.~E. Dullerud, ``A unified framework for
  lp analysis and synthesis of linear switched systems,'' in \emph{2016
  American Control Conference (ACC)}.\hskip 1em plus 0.5em minus 0.4em\relax
  IEEE, 2016, pp. 715--720.

\bibitem{PeBaRTOD}
M.~Petreczky, L.~Bako, and J.~H. Van~Schuppen, ``Realization theory of
  discrete-time linear switched systems,'' \emph{Automatica}, vol.~49, no.~11,
  pp. 3337--3344, 2013.

\bibitem{PeScRTOD}
M.~Petreczky and J.~H. van Schuppen, ``Realization theory of discrete-time
  linear hybrid system,'' \emph{IFAC Proceedings Volumes}, vol.~42, no.~10, pp.
  593--598, 2009.

\bibitem{PhEsSODT}
M.~Philippe, R.~Essick, G.~Dullerud, and R.~M. Jungers, ``Stability of
  discrete-time switching systems with constrained switching sequences,''
  \emph{Automatica}, vol.~72, pp. 242--250, 2016.

\bibitem{PhMiDTBA}
M.~Philippe, G.~Millerioux, and R.~M. Jungers, ``Deciding the boundedness and
  dead-beat stability of constrained switching systems,'' \emph{Nonlinear
  Analysis: Hybrid Systems}, 2016.

\bibitem{PuZhASOT}
V.~Putta, G.~Zhu, J.~Shen, and J.~Hu, ``A study of the generalized
  input-to-state l2-gain of discrete-time switched linear systems,'' in
  \emph{50th IEEE Conference on Decision and Control and European Control
  Conference}, 2011, pp. 435--440.

\bibitem{ShWiAPSM}
R.~Shorten, F.~Wirth, and D.~Leith, ``A positive systems model of tcp-like
  congestion control: asymptotic results,'' \emph{IEEE/ACM Transactions on
  Networking}, vol.~14, no.~3, pp. 616--629, 2006.

\end{thebibliography}

\end{document}